# О вычислении топологической степени оператора Пуанкаре в сингулярных случаях[1]

И. С. Мартынова

Воронежский Государственный Университет

**Аннотация:** В статье развивается метод оценки топологической степени отображения Пуанкаре математической модели мелкой лагуны, находящейся под T-периодическим воздействием. С помощью этого метода мы получаем условия для параметров модели, гарантирующие существование T-периодического решения в определенной области. Сложность задачи связана с тем, что непосредственное применение принципа невозвращаемости М.А. Красносельского-А.И. Перова приводит к сингулярному векторному полю.

# A geometrical method to study plankton oscillations in the mathematical model of a shallow lagoon under periodic climate change

S. Martynova
Voronezh State University

**Annotation:** This paper develops a technique to evaluate the topological degree of the Poincare map associated to the mathematical model of a shallow lagoon subjected to a T-periodic forcing. By means of this technique we derive conditions for the parameters of the model guaranteeing the existence of a T-periodic solution in a certain region.

## 1. Введение.

В статье изучается существование 12-периодических решений системы, описывающей математическую модель мелкой лагуны при периодическом изменении климата:

$$x'_1 = \lambda_1 s(t) x_1 - \lambda_2 x_1^2 - \lambda_3 x_2 \left( \frac{x_1}{k_P + x_1} \right), \qquad (1)$$

$$x'_2 = \lambda_4 x_2 \left( \frac{x_1}{k_P + x_1} \right) - \lambda_5 x_2, \qquad (2)$$

где дополнительное внешнее воздействие имеет вид:

$$s(t) = M + N \sin\left( \frac{2\pi}{12} t + 1 \right). \qquad (3)$$

Рассматриваемая задача обусловлена рядом численных экспериментов на тему периодических решений системы (1)-(2), опубликованных в статье [1]. Целью данной статьи является дать некоторое теоретическое обоснование полученных в [1] численных результатов. Статья организована следующим образом. В следующей части приводится основной результат о топологической степени оператора Пуанкаре системы (1)-(2). В заключительной части основной результат применяется для доказательства существования 12-периодического решения этой системы с некоторыми параметрами, обусловленными динамикой модели мелкой лагуны.

Для дальнейшего рассмотрения введем следующие константы:

$$u_{N,1} = \frac{\lambda_1}{\lambda_2}(M - N),$$

$$p_{N,1} = \frac{\lambda_1}{\lambda_2}(M + N),$$

$$r_1 = \frac{\lambda_5 k}{\lambda_4 - \lambda_5},$$

$$r_{N,2} = (k + r_1)(\lambda_1 (M + N) - \lambda_2 r_1).$$

## 2. Геометрический метод для получения результатов о динамике модели лагуны.

Следующая теорема дает достаточные условия существования в системе 12-периодических решений, лежащих строго в первой четверти.

**Теорема.** Предположим, что $\lambda_i > 0,\ i \in \overline{1,5},\ k > 0,\ M > 0,\ \lambda_4 > \lambda_5$ и зафиксируем такое $N_0 \geq 0$, что

$$\frac{\lambda_1}{\lambda_2}(M - N_0) > \Gamma_1. \qquad (4)$$

Тогда существует $\eta > 0$ такое, что при каждом $N \in [0, N_0]$ система (1)-(3) имеет 12-периодическое решение $x_N$, удовлетворяющее соотношениям

$$\eta \le x_{N,1} \le p_1 + \eta,$$

$$\eta \le x_{N,2} \le \max\{1, \Gamma_{N,2}\} e^{12\lambda_4 \frac{p_{N,1}}{k+p_{N,1}} + 12\lambda_5} + \eta.$$

**Доказательство.** Рассмотрим возмущенную систему

$$\begin{cases} x'_1 = \lambda_1\left(M + N\sin\left(\dfrac{2\pi}{12}t + 1\right)\right)x_1 - \lambda_2 x_1^2 - \\ -\lambda_3 x_2 \dfrac{x_1}{k + x_1} + \varepsilon := P_{\varepsilon,N}(t, x_1, x_2), \\ x'_2 = \lambda_4 x_2 \dfrac{x_1}{k + x_1} - \lambda_5 x_2 x_2^{\varepsilon} := Q_{\varepsilon,N}(x_1, x_2), \end{cases} \qquad (5)$$

где $\varepsilon > 0$ и обозначим через $\Omega_{\varepsilon,N}(T)$ оператор Пуанкаре-Андронова системы (5). Положим

$$f_{\varepsilon,N,1}(x_1) = \left((\lambda_1(M - N) - \lambda_2 x_1)(k + x_1) + \varepsilon + \frac{\varepsilon k}{x_1}\right),$$

$$i \in \overline{1,2}, N \in [0, N_0],$$

$$f_{\varepsilon,N,2}(x_1) = \left((\lambda_1(M - N) - \lambda_2 x_1)(k + x_1) + \varepsilon + \frac{\varepsilon k}{x_1}\right),$$

$$i \in \overline{1,2}, N \in [0, N_0],$$

$$g_{\varepsilon}(x_1) = \left(\frac{\lambda_4}{\lambda_5} \cdot \frac{x_1}{k + x_1}\right)^{\frac{1}{\varepsilon}}.$$

Можно проверить, что для $x \in R_+^2$ и $t \in [0, T]$ выполнены следующие соотношения:

$$P_{\varepsilon,N}(t, x_1, x_2) \le 0$$

для любого $x_2 \ge f_{\varepsilon,N,2}, N \in [0, N_0],$ (6)

$$P_{\varepsilon,N}(t, x_1, x_2) \ge 0$$

для любого $x_2 \le f_{\varepsilon,N,1}, N \in [0, N_0],$ (7)

$$G_{\varepsilon}(x_1, x_2) \le 0$$

для любого $x_2 \ge g_{\varepsilon}(x_1),$ (8)

$$G_{\varepsilon}(x_1, x_2) \ge 0$$

для любого $x_2 \le g_{\varepsilon}(x_1).$ (9)

Пусть $\eta > 0$ произвольная фиксированная константа, рассмотрим $\varepsilon \in (0, \eta)$. Заметим, что существует единственное положительное число $p_{\varepsilon,N,1}$ такое, что $f_{\varepsilon,N,2}(p_{\varepsilon,N,1}) = 0$. Следовательно, мы имеем

$p_{\varepsilon,N,1} \to p_{N,1}$ при $\varepsilon \to 0$, (10)

равномерно по $N \in [0, N_0]$.

Пусть $p_{\varepsilon,N,2} \in (g_\varepsilon(p_{\varepsilon,N,1}), \infty)$ и обозначим через $e_{\varepsilon,N,1}$ такую точку, что $f_{\varepsilon,N,1}(e_{\varepsilon,N,1}) > p_{\varepsilon,N,2}$ и $e_{\varepsilon,N,1} \in (0,\varepsilon)$. Пусть $e_{\varepsilon,N,1} \in (0, g(e_{\varepsilon,N,1}))$. Заметим, что $e_{\varepsilon,N,1} < p_{\varepsilon,N,1}, e_{\varepsilon,N,2} < p_{\varepsilon,N,2}$ и

$$e_{\varepsilon,N,1} \to 0, e_{\varepsilon,N,2} \to 0 \text{ как } \varepsilon \to 0 \qquad (11)$$

равномерно по $N \in [0, N_0]$.

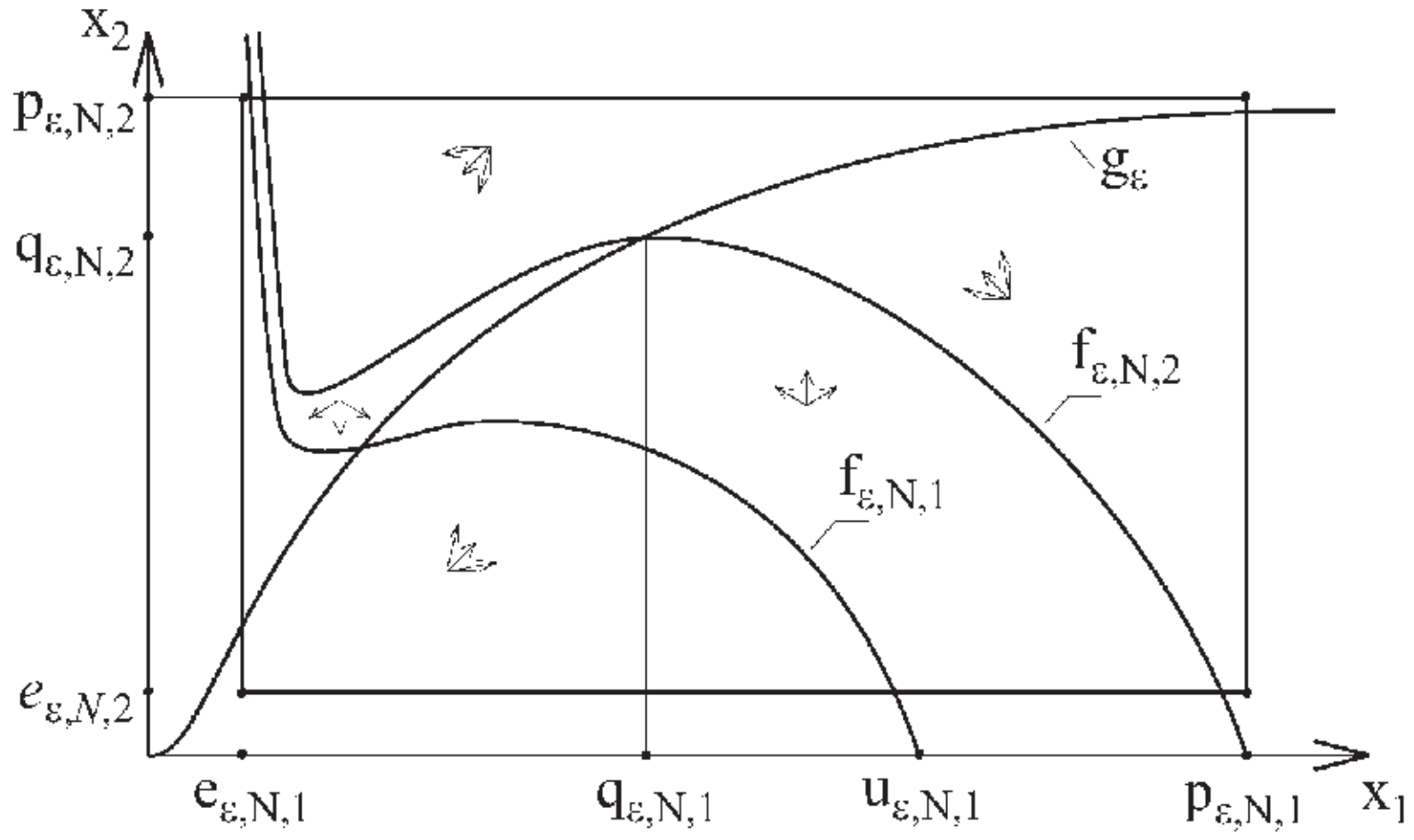

Рис. 1: Схематичный график векторного поля $(P_\varepsilon(t, x_1, x_2), Q_\varepsilon(x_1, x_2))$ при $t \in [0, T]$, где различные направления вектора в фиксированной точке соответствуют различным значениям $t \in [0, T]$.

По построению (см. рис. 1) граница прямоугольного множества $R_{\varepsilon,N} = (e_{\varepsilon,N,1}, e_{\varepsilon,N,2}, p_{\varepsilon,N,1}, p_{\varepsilon,N,2})$ содержит точки T - невозвращаемых решений системы (5) и, следовательно (см. [2], лемма 6.1),

$$d(I - \Omega_{\varepsilon,N}, R_{\varepsilon,N}) = d\left(-\begin{pmatrix} P_{\varepsilon,N}(0,\cdot) \\ G_{\varepsilon,N} \end{pmatrix}, R_{\varepsilon,N}\right).$$

С другой стороны, из соотношений (6)-(9) следует, что поле

$$-\begin{pmatrix} P_{\varepsilon,N}(0,\cdot) \\ G_{\varepsilon,N} \end{pmatrix}, (R_{\varepsilon,N})$$

на границе множества $R_{\varepsilon,N}$ направлено во внутрь (см. рис. 2) и, таким образом,

$$d(I - \Omega_{\varepsilon,N}, R_{\varepsilon,N}) = d\left(-\begin{pmatrix} P_{\varepsilon,N}(0,\cdot) \\ G_{\varepsilon,N} \end{pmatrix}, R_{\varepsilon,N}\right) = 1 \qquad (12)$$

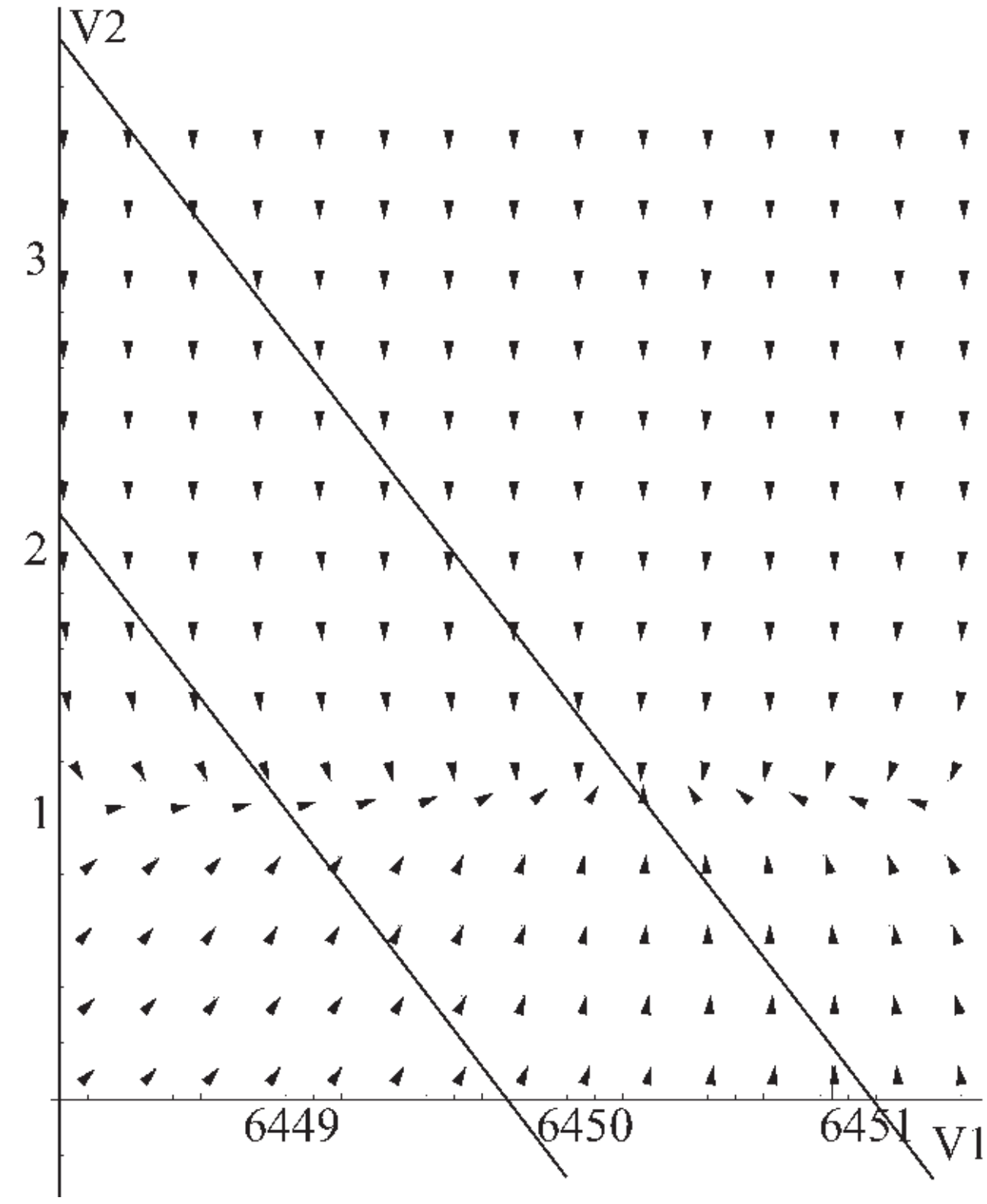

Рис. 2: Схематичный график векторного поля $(P_{\varepsilon,N}(0,\cdot),G_{\varepsilon,N})$ при $\varepsilon=0,3$, значениях $\lambda_1=1.29, \lambda_2=0.0002, \lambda_3=0.93, \lambda_4=6.9,$ $\lambda_5=5.8, k=100$.

Заметим, что мы можем выбрать достаточно малое $\varepsilon_0\in(0,\eta)$ такое, что для любого $\varepsilon\in(0,\varepsilon_0)$ будет существовать единственная константа $q_{\varepsilon,N,1}>0$ такая, что $f_{\varepsilon,N,2}(q_{\varepsilon,N,1})=g_\varepsilon(q_{\varepsilon,N,1})$. Пусть $q_{\varepsilon,N,2}=f_{\varepsilon,N,2}(q_{\varepsilon,N,1})$. Так как

$$f_{\varepsilon,N,2}(x_1)\to f_{0,N,2}(x_1) \text{ при } \varepsilon\to 0, x_1>0, \qquad (13)$$

$$g_\varepsilon(x_1)\to 0 \text{ при } \varepsilon\to 0, x_1\in(0,r_1), \qquad (14)$$

$$g_\varepsilon(x_1)\to\infty \text{ при } \varepsilon\to 0, x_1>r_1, \qquad (15)$$

и $f_{0,N,2}(r_1)=r_{N,2}>0$, то

$$q_{\varepsilon,N,1}\to r_1 \text{ при } \varepsilon\to 0.$$

Следовательно,

$$q_{\varepsilon,N,2}\to r_2 \text{ при } \varepsilon\to 0. \qquad (16)$$

Обозначим через $x_{\varepsilon,N}$ 12-периодическое решение системы (5) такое, что $x_{\varepsilon,N}(0)\in R_{\varepsilon,N}$. Покажем, что

$$x_{\varepsilon,N,2}(0)\le \max\{1,q_{\varepsilon,N,2}\}e^{12\lambda_4\frac{p_{\varepsilon,N,1}}{k+p_{\varepsilon,N,1}}+12\lambda_5} \qquad (17)$$

для любого $\varepsilon\in(0,\varepsilon_0)$ и $N\in[0,N_0]$.

Предположим противное, что существует $\varepsilon\in(0,\varepsilon_0)$ и $N\in[0,N_0]$ такие, что

$$x_{\varepsilon,N,2}(0) > \max\{1, q_{\varepsilon,N,2}\} e^{12\lambda_4 \frac{p_{\varepsilon,N,1}}{k+p_{\varepsilon,N,1}} + 12\lambda_5}. \qquad (18)$$

Продолжим функцию $x_{\varepsilon,N}$ на интервал [-12,0], периодически. Из геометрических соображений существует $s \in [-12{,}0)$ такое, что

$$\begin{cases} x_{\varepsilon,N,2}(s) = q_{\varepsilon,N,2}, \\ x_{\varepsilon,N,2}(t) > q_{\varepsilon,N,2},\ t \in (s,0]. \end{cases} \qquad (19)$$

Действительно, если $x_{\varepsilon,N,2}(t) = g_{\varepsilon,N,2}$ для всех $t \in [0{,}12]$, то условия (6)-(9) означают, что неравенства $x'_{\varepsilon,N,1} > 0$ и $x'_{\varepsilon,N,2} > 0$ не могут быть выполнены одновременно при $t \in [0{,}12]$. Но это противоречит тому, что система (5) не имеет постоянных решений и производная $x'_{\varepsilon,N}$ непрерывно-дифференцируемой 12-периодической непостоянной функции $x_{\varepsilon,N}$, действующей в $R^2$, должна пройти все направления при t меняющемся от 0 до 12. Из (18) имеем, что $x_{\varepsilon,N,2}(0) > 1$ и, учитывая (19), выводим (см. [2], § 1.4) существование $\tau \in [s,0)$ такого, что

$$\begin{cases} x_{\varepsilon,N,2}(\tau) = \max\{1, q_{\varepsilon,N,2}\}, \\ x_{\varepsilon,N,2}(t) > \max\{1, q_{\varepsilon,N,2}\},\ t \in (\tau,0]. \end{cases} \qquad (20)$$

Так как

$$\begin{aligned} &\left(\lambda_4 x_2 \frac{x_1}{k+x_1} - \lambda_5 x_2\right) - \\ &-\left(\lambda_4 x_2 \frac{x_1}{k+x_1} - \lambda_5 x_2 x_2^{\varepsilon}\right) = \qquad (21) \\ &= \lambda_5 x_2 x_2^{\varepsilon} - \lambda_5 x_2 = \\ &= \lambda_5 x_2 (x_2^{\varepsilon} - 1) > 0,\ x_2 > 1 \end{aligned}$$

тогда, используя (20) для решения $y_{\varepsilon,N,\delta}$ задачи Коши

$$\begin{cases} y' = \left(\lambda_4 \frac{x_{\varepsilon,N,1}(t)}{k + x_{\varepsilon,N,1}(t)} - \lambda_5\right) y, \\ y(\tau) = x_{\varepsilon,N,2}(\tau) + \delta, \end{cases} \qquad (22)$$

где $\delta > 0$, имеем неравенство

$$x_{\varepsilon,N,2} < y_{\varepsilon,N,\delta},\ t \in [\tau,0].$$

Переходя к пределу в последнем неравенстве, получаем

$$x_{\varepsilon,N,2} < y_{\varepsilon,N,0},\ t \in [\tau,0],$$

где $y_{\varepsilon,N,0}$ есть решение (22) с $\delta = 0$. Таким образом, имеем оценку

$$x_{\varepsilon,N,2}(0) \le y_{\varepsilon,N}(0) =$$
$$= x_{\varepsilon,N,2}(\tau) e^{\int_\tau^0 \left(\lambda_4 \frac{x_{\varepsilon,N,1}(t)}{k + x_{\varepsilon,N,1}} - \lambda_5\right) d\tau} \le$$
$$\le \max\{1, q_{\varepsilon,N,2}\} e^{12\lambda_4 \frac{p_{\varepsilon,N,1}}{k+p_{\varepsilon,N,1}} + 12\lambda_5},$$

которая удовлетворяет (18) и, следовательно, (17) выполнено. Используя отношения (10) и (16), из (17) получаем, что

$$x_{\varepsilon,N,2}(0) < \max\{1, r_{N,2}\} e^{12\lambda_4 \frac{p_{\varepsilon,N,1}}{k+p_{\varepsilon,N,1}} + 12\lambda_5} + \eta, \quad (23)$$
$$\varepsilon \in (0, \varepsilon_0), N \in [0, N_0].$$

Покажем, что $\eta > 0$ может быть уменьшено настолько, что

$$x_{\varepsilon,N,1}(0) \ge \eta \text{ для всех } \varepsilon \in (0, \varepsilon_0), N \in [0, N_0]$$

и

$$x_{\varepsilon,N,2}(0) \ge \eta \text{ для всех } \varepsilon \in (0, \varepsilon_0), N \in [0, N_0].$$

Предположим противное, то есть, что существуют последовательности

$$\varepsilon_k \to 0, N_k \to N_* \in [0, N_0] \text{ при } k \to \infty$$

такие, что

$$x_{\varepsilon,N_k,1}(0) \to 0 \text{ при } k \to \infty \quad (24)$$

или

$$x_{\varepsilon,N_k,2}(0) \to 0 \text{ при } k \to \infty. \quad (25)$$

Так как $x_{\varepsilon_k,N_k}$ удовлетворяет (5), последовательность $\{x_{\varepsilon_k,N_k}\}_{k=1}^{\infty}$ - компактна и без ограничения общности можем рассматривать $x_{\varepsilon_k,N_k} \to x_{0,N_*}$ при $k \to \infty$.

Из (6) и (7) (см. рис. 1) имеем существование некоторого $t_n \in [0,12]$, такого, что

$$x_{\varepsilon_k,N_k,2}(t_n) \in$$
$$\in [f_{\varepsilon_k,N_k,1}(x_{\varepsilon_k,N_k,1}(t_k)), f_{\varepsilon_k,N_k,2}(x_{\varepsilon_k,N_k,1}(t_k))].$$

На самом деле, одно из следующих отношений выполнено

$$x'_{\varepsilon_k,N_k,1}(t) > 0, t \in [0,12]$$

или

$$x'_{\varepsilon_k,N_k,1}(t) < 0, t \in [0,12].$$

Но, как уже было замечено, вектор $x'_{\varepsilon,N}(t)$ должен пройти все направления при t меняющемся от 0 до 12. Аналогично, существует $s_n \in [0,12]$ такое, что

$$x_{\varepsilon_k,N_k,2}(s_n) = g_{\varepsilon_n}(x_{\varepsilon_k,N_k,1}(s_n)).$$

Следовательно, по условиям (13), (15) существуют $t_0, s_0 \in [0,12]$ такие, что

$$x_{0,N_*,2}(t_0) \in \in [f_{0,N_*,1}(x_{0,N_*,1}(t_0)), f_{0,N_*,2}(x_{0,N_*,1}(t_0))], \quad (26)$$

$$x_{0,N_*,1}(s_0) \le r_{1.} \quad (27)$$

Из соотношений (24)-(25) следует, что либо

$$x_{0,N_*,1}(0) = 0, \quad (28)$$

либо

$$x_{0,N_*,2}(0) = 0, \quad (29)$$

Покажем, что оба этих равенства не верны.
Во-первых, равенство (28) влечет

$$x_{0,N_*,1}(t) = 0, t \in [0,12]$$

и из (26) имеем, что $x_{0,N_*,2}(t) \ne 0$. Но можно заметить, что система (1) не имеет ненулевых 12-периодических решений вида $x_0(t) = (0, x_{02}(t))$. Аналогично, если имеет место равенство (29), то

$$x_{0,N_*,2}(t) = 0, t \in [0,12]. \quad (30)$$

Обозначим через $u_{\varepsilon,N,1}$ такую точку, при которой $f_{\varepsilon,N,1}(u_{\varepsilon,N,1}) = 0$.
Заметим, что

$$u_{\varepsilon_k,N_k,2} \to u_{1,N} \text{ при } k \to \infty, \quad (31)$$

из (30) и (26) получаем

$$x_{0,N_*,1}(t_0) \ge u_{1,N}. \quad (32)$$

Также из (31) и (7) имеем, что

$$P_{0,N_*}(t, x_1, 0) \ge 0, x_1 \in [0, u_{1,N_*}], t \in [0,12]. \quad (33)$$

Из периодичности функции $x_{0,N_*,1}$, соотношений (32) и (33) заключаем

$$x_{0,N_*,1}(t) \ge u_{1,N_*}, t \in [0,12]. \quad (34)$$

Оценки (27) и (34) противоречат условию (4) и, таким образом, (29) не имеет места.

Доказательство завершено.

## 3. Приложение полученного результата.

Используя следующие значения параметров системы (1)-(2), (см. [1], таблицы 3, 4):

$$\lambda_1 = 1.29, \lambda_2 = 0.0002, \lambda_3 = 0.93,$$
$$\lambda_4 = 6.9, \lambda_5 = 5.8, k = 100,$$

получаем следующее следствие доказанной теоремы.

**Следствие.** Пусть параметры системы (1)-(2) заданы значениями

$$\lambda_1 = 1.29, \lambda_2 = 0.0002, \lambda_3 = 0.93,$$
$$\lambda_4 = 6.9, \lambda_5 = 5.8, k = 100,$$

и

$$M > \frac{0,116}{1,419}.$$

Выберем какое-нибудь $0 < N_0 < M - \frac{0,116}{1,419}$. Тогда существует $\eta > 0$ такое, что при любом $0 \le N \le N_0$ система (1)-(2) допускает по крайней мере одно 12-периодическое решение $x_N$ с начальным условием, удовлетворяющим соотношениям

$$\eta \le x_{N,1} \le \frac{1,29}{0,0002}(M+N) + \eta,$$
$$\eta \le x_{N,2} \le$$
$$\le \max\left\{1, \left(100 + \frac{580}{1,1}\right) \cdot \left(1,29(M+N) - \frac{0,116}{1,1}\right)\right\} \cdot$$
$$\cdot e^{\frac{534060(M+N)}{100+6450(M+N)} + 69,6} + \eta.$$

**Доказательство.** Достаточно заметить, что условие (4) в случае указанных значениях параметров принимает вид

$$M - N_0 > \frac{0,116}{1,419}.$$

Литература

[1] A. Garulli, C. Mocenni, A. Tesi, A. Vicino. Integrating identification and qualitative analysis for the dynamic model of a lagoon, Int. J. Bif. Chaos, Vol.13, No.2 (2003) 357-374.

[2] М. А. Красносельский, Оператор сдвига по траекториям дифференциальных уравнений, 1968.